\documentclass{article}
\usepackage[utf8]{inputenc}
\usepackage{tabularx}
\usepackage{amsmath}
\usepackage{amssymb}
\usepackage{amsthm}
\usepackage{geometry}
\usepackage{enumitem}
\usepackage[T1]{fontenc}
\usepackage[dvipsnames]{xcolor}
\title{On equations of fifth degree  }
\author{by Mr. Kronecker, Translation by Yonathan Stone\thanks{The translator would like to thank Jesse Wolfson for suggesting this and for providing helpful feedback on earlier drafts. In addition, this translation was supported in part by NSF Grant DMS-1944862.}}
\date{from the monthly report of the Academy of Sciences at Berlin, 1861, translation from 2020}

\begin{document}
\maketitle
I have very recently completed a problem that has occupied me from time to time these last five years.  I continued to take it up time and time again because its resolution would prove decisive for the future direction of my algebraic investigations.  That is to say that very early on in my work concerning the algebraic solutions of equations I arrived at the insight that the problem can be viewed more generally from two different sides.  On the one hand one can assume more general rational functions of the roots, which I will call \textit{Affectfunctionen}, in place of the coefficients of the equations, that is symmetric functions of the same roots.  On the other hand one can consider more general algebraic functions in place of the usual root symbols, that is instead of those function symbols that are purely defined by the equations, and use these as auxiliary functions in finding a solution. The first side of this more general viewpoint has already been adopted in older algebraic work, even if this was not clearly stated or indicated.  However, the second type of generalization of this problem has only recently been taken up simultaneously by Mr. Hermite as well as myself.  Even so, the way in which I chose to pursue this was markedly different from that used by Mr. Hermite.  In particular, a postulate I made regarding the method of solution complicated answering the question for equations of fifth degree, but this in turn led to further interesting investigations due to be grounded in nature of the problem.
\\
``If an equation is algebraically solvable, one can always give the root a form that allows any algebraic function it is made up of to be expressed using rational functions of the roots of said equation.'' This is the statement of Abel's theorem, which contains an extremely important property of the ordinary root expressions.  This property is also what must continue to hold for the more general expressions of roots belonging to equations that are not solvable in the ordinary sense.  And thanks to this postulate we must choose the algebraic function symbols that help solving equations in the broad sense.
\\
In the course of investigating equations of fifth degree, the aforementioned postulate, which I will justify in a more comprehensive communication, led me to looking for rational functions of the five roots with the property that the distinct values under permutation of the roots possessed as many identical relations as possible.  For reasons that are abundantly clear I decided to only consider those permutations which left the square root of the discriminant of the equation of fifth degree unchanged. By doing so I in fact found twelve-valued rational functions of the five roots with the property that each of the twelve values differs from another value by only a sign and that the six distinct absolute values are linked to one another by three linear relations.  For that reason the square of such a function is the root of an equation of sixth degree, and the coefficients thereof are rationally constructed from those of the original equation and the square root of the discriminant and only dependent on three such rational expressions.  In discovering these kinds of functions I at first nearly succeeding at using the modular equation of fifth degree to solve the equation of fifth degree with almost no computations.  For this reason I include two such remarkable functions in a letter to Mr. Hermite, which has been reprinted in the 1858 \textit{comptes rendus} of the Paris Academy.  Secondly this provided the possibility of solving general equations of fifth degree in one of the ways corresponding to the postulate above, admittedly only by using algebraic functions of two variables.  To discuss this important point somewhat closer, I will set:
\[f(x_0,x_1,x_2,x_3,x_4) = \sum\sum \sin\left(\frac{2n\pi}{5}\cdot x_m x^2_{m+n} x^2_{m + 2n},\right)\]
where $x_0,x_1,x_2,x_3,x_4$ are roots of an equation $X = 0$ of fifth degree, the summations run from $m = 0,1,2,3,4$ and $n = 1,2,3,4$, and larger indices are reduced to the smallest remainders modulo 5.  It thereupon follows that $f(x_0,x_1,x_2,x_3,x_4)$ satisfies an equation of twelfth degree:
\[(f^2 + a)^6 + 4a(f^2 + a)^5 + 10b(f^2 + a)^3 + 4c(f^2 + a) - 4ac + 5b^2 = 0, \tag{I.}\]
where $a,b,c$ are two-valued polynomial functions of the five roots $x$.  But there also exist countless other rational functions of the roots $x$ besides the function $f$ which share the property of satisfying twelfth degree equations of the indicated form\footnote{Regarding this one may also look at the work of Mr. Brioschi in his address: \textit{``Sul metodo di Kronecker per la risolusione delle equasioni di quinto grado''} [On Kronecker's method for solving equations of fifth degree] (read on the 25. Novemeber 1858 at the Lombardy Institute), where for a particular function $f$ the complete expression of the coefficients $a,b,c$ is initially also given by the invariants of the equation of fifth degree.}, and among the rational functions of this kind there are in turn those for which the expressions corresponding to the values $a,b,c$ depend only on two rational two-valued functions: $\varphi(x_0,x_1,x_2,x_3,x_4), \psi(x_0,x_1,x_2,x_3,x_4).$  Such a specialized function $f$ is thus an implicitly given algebraic function in $\varphi$ and $\psi$ may therefore be designated by $W(\varphi,\psi)$.  Since this function $f$ is cyclic and thus allows the equation $X = 0$ to be solvable, we have that the roots of the general equation of fifth degree admit an explicit representation using square roots, fifth roots, and the function symbol $W$, and in particular in such a way that fully satisfies the postulate above.  All of this more or less became clear during the discovery of the functions $f$ as an immediate consequence.  But it remained to determine whether this question was completed, that is, it was not yet clear whether the two variable algebraic function $W$ could be reduced to one of a \textit{single} variable.  It has been known for a while that such a reduction is possible, provided one drops the oft-mentioned postulate and has been recently presented in the letter to Mr. Hermite.  From further occupation with these objects I have obtained more specifically relevant results.  But only recently did I suceed to answer the main question and ascertained that the reduction of the algebraic function $W$ to a function of \textit{one} variable and therefore the solution of general equation of fifth degree using algebraic functions of one variable is indeed impossible, provided that Abel's theorem introduced above continues to hold, which is in fact the case for the solution of equations through root symbols.  It seems to me that this result constitutes an impressive extension of  Abel's proof of the impossibility of algebraically solving equations of higher degree;  and it simultaneously contains the conclusion of the fifth degree solvability problem in its more general setup.  This is a conclusion for which  I cannot view my results as finished and suitable for publishing prior to its achievement.
\\
That in the subject of algebra it becomes neccessary to introduce two variable functions, and then in a certain sense reduce these back to functions of a \textit{single} variable, should not present a cause for concern.  The reader should recall that in Analysis the fourfold periodic functions of two variables were introduced by Jacobi and are definitely still used, in spite of the fact that in the 30\textsuperscript{th} volume of this journal he himself showed that these may be constructed algebraically from functions of \textit{one} variable.  Without delving further into this analogy I would like to return to equations of fifth degree and attach a few more remarks to the method of solution alluded to above.
\\
As before, letting $f(x_0,x_1,x_2,x_3,x_4)$ denote a function satisfying an equation of the form (I.), and setting
\[f_k = f(x_k,x_{k+3},x_{k+4},x_{k+1},x_{k+2}),\]
it follows that $\pm f, \pm f_0, \pm f_1, \pm f_2, \pm f_3, \pm f_4$ are the roots of the aforementioned equation.  In turn there exist rational functions of the six distinct quantities $f$ which themselves are roots of equation of fifth degree whose coefficient are rationally constructed from $a,b$ and $c$.   Letting $\Phi$ denote one of the these functions of the quantities $f$, it follows that $\Phi$ is simultaneously a rational function of the roots $x$ themselves.  Moreover, when presented as a polynomial in \textit{one} of the roots $x$, its coefficients only contain rational expressions made up of the coefficients of the equation $X = 0$ and the square root of its discriminant.  In addition it can be shown without additional computation that the product $(f-f_0)(f_1-f_4)(f_2-f_3)$, which is included among the functions denoted by $\Phi$, satisfies an equation of fifth degree, for which both the second and the fourth coefficient are equal to zero.  This result, which is of a certain formal interest, can otherwise be found in the notes with which Mr. Brioschi accompanied the announcement of Hermite's solution of the equations of fifth degree.  Additionally, in a more recent communication by letter to Mr. Borchardt\footnote{pg. 304}, Mr. Hermite indicated a special function of the roots of an equation of fifth degree, also belonging to the functions $\Phi$, which bears special interest due to both its simplicity as well as its relation to the invariants of the equation.  The essence of this matter, which arises from the simplest investigations of the functions $f$, can be summarized in the following way:
\newpage
\begin{center}
    Among the ten-valued rational functions of the five quantities $x_0,x_1,x_2,x_3,x_4$, which take on only five values under all cyclic permutations of exactly 3 of the these quantities, there are those for which the symmetric functions of the five values only depend on two functions of the quantities $x$.  Morever, among these there are especially those for which the sum of the five values themselves as well as the sum of third powers identically vanish.
\end{center}
Algebraic functions that essentially depend on two variable are defined through the hereafter appearing \textit{special} equations of fifth degree, for instance the equations of the form $z^5 + pz^3 + qz + r$.  That is, these introduce algebraic functions that may be used for the solution of the \textit{general} equation of fifth degree, much in the same way as the function $W$ above.  But with respect to certain more general solvability problems the twelfth order equation introduced above is preferred as an auxiliary equation.  And a more accurate discussion of this allows one to recognize its many remarkable properties and thus simultaneously the different forms one may give the roots of equations of fifth degree using the symbol $W$.
\\
Berlin, June 1861
\end{document}